\documentclass[12pt,a4paper]{article}

\usepackage{amsfonts}
\usepackage{mathrsfs}
\usepackage{indentfirst, latexsym, bm}
\usepackage{amsfonts,amssymb,amsthm}
\usepackage{amsbsy}
\usepackage{latexsym}
\usepackage[mathscr]{eucal}
\usepackage{amsmath,amscd}  %,mathrsfs'''''''''

\newtheorem{theorem}{\hspace{1.3em}Theorem}[section]
\newtheorem{lemma}{\hspace{1.3em}Lemma}[section]

\newtheorem{corollary}{\hspace{1.3em}Corollary}[section]
\newtheorem{remark}{\hspace{1.3em}Remark}[section]
\newtheorem{example}{\hspace{1.3em}Example}[section]

\begin{document}

\title{Ordered $k$-flaw Preference Sets}

\author{Po-Yi Huang$^{a,}$\thanks{Email address of the corresponding author:
pyhuang@mail.ncku.edu.tw
 Partially supported by NSC 96-2115-M-006-012
}
 \and Jun Ma$^{b}$
 \and Yeong-Nan Yeh$^{b,}$\thanks{Partially supported by NSC 96-2115-M-006-012}}

\date{}
\maketitle \vspace{-1cm} \begin{center} \footnotesize
 $^{a}$ Department of Mathematics,
National Cheng Kung University, Tainan, Taiwan
\\
$^{b}$ Institute of Mathematics, Academia Sinica, Taipei, Taiwan
\end{center}

\thispagestyle{empty}\vspace*{.4cm}

\begin{abstract}
In this paper, we focus on ordered $k$-flaw preference sets. Let
$\mathcal{OP}_{n,\geq k}$ denote the set of ordered preference sets
of length $n$ with at least $k$ flaws and
$\mathcal{S}_{n,k}=\{(x_1,\cdots,x_{n-k})\mid x_1+x_2+\cdots
+x_{n-k}=n+k, x_i\in\mathbb{N}\}$. We obtain a bijection from the
sets $\mathcal{OP}_{n,\geq k}$ to $\mathcal{S}_{n,k}$. Let
$\mathcal{OP}_{n,k}$ denote the set of ordered preference sets of
length $n$ with exactly $k$ flaws. An $(n,k)$-\emph{flaw path} is a
lattice path starting at $(0,0)$ and ending at $(2n,0)$ with only
two kinds of steps---rise step: $U=(1,1)$ and fall step: $D=(1,-1)$
lying on the line $y = -k$ and touching this line. Let
$\mathcal{D}_{n,k}$ denote the set of $(n, k)$-flaw paths. Also we
establish a bijection between the sets $\mathcal{OP}_{n,k}$ and
$\mathcal{D}_{n,k}$. Let $op_{n,\geq k,\leq l}^m$ $(op_{n, k,
=l}^m)$ denote the number of preference sets
$\alpha=(a_1,\cdots,a_n)$ with at least $k$ (exact) flaws and
leading term $m$ satisfying $a_i\leq l$ for any $i$ $(\max\{a_i\mid
1\leq i\leq n\}=l)$, respectively. With the benefit of these
bijections, we obtain the explicit formulas for $op_{n,\geq k,\leq
l}^m$. Furthermore, we give the explicit formulas for $op_{n, k,
=l}^m$. We derive some recurrence relations of the sequence formed
by ordered $k$-flaw preference sets of length $n$ with leading term
$m$. Using these recurrence relations, we obtain the generating
functions of some corresponding $k$-flaw preference sets.
\end{abstract}

\noindent {\bf Keyword: Dyck path; flaw; Parking function;
Bijection; Catalan number}

\newpage

\section{Introduction}
Throughout the paper, we let $[m,n]:=\{m,m+1,\cdots,n\}$,
$[n]:=\{1,2,\cdots,n\}$ and $c_n$ denote the $n$-th Catalan number
for any $n\geq 0$, i.e. $c_n=\frac{1}{n+1}{2n\choose n}$. It is well
known that the sequence $c_n$ has the generating function
$C(x):=\sum\limits_{n\geq 0}c_nx^n$ satisfying the functional
equation $C(x)=1+zC(x)^2$ and $C(x)=\frac{1-\sqrt{1-4x}}{2x}$
explicitly.

Suppose that $n$ cars will be parked in $n$ parking spaces which are
arranged in a line and numbered 1 to n from left to right. Each car
has initial parking preference $a_i$. If space $a_i$ is occupied,
the car moves to the first unoccupied space to the right. We call
$(a_1,\cdots,a_n)$ {\it preference set}. Clearly, the number of the
preference sets is $n^n$. If a preference set $(a_1,\cdots,a_n)$
satisfies $a_i\leq a_{i+1}$ for $1\leq i\leq n-1$, then we say that
this preference set is {\it ordered}. If all the cars can find the
parking spaces, then we say the preference set is a {\it parking
function}. If there are exactly $k$ cars which can't be parked, then
the preference set is called a {\it $k$-flaw preference set}.

Parking functions were introduced by J. Riordan in \cite{R}. He
derived that the number of parking functions of length $n$ is
$(n+1)^{n-1}$, which coincides with the number of labeled trees on
$n+1$ vertices by Cayley's formula. Several bijections between the
two sets are known (e.g., see \cite{FR,R,SMP}). Furthermore, Riordan
concluded that the number of ordered parking functions is
$\frac{1}{n+1}{2n\choose{n}}$, which is also equals the number of
Dyck path of semilength $n$. Parking functions have been found in
connection to many other combinatorial structures such as acyclic
mappings, polytopes, non-crossing partitions, non-nesting
partitions,  hyperplane arrangements, etc. Refer to
\cite{F,FR,GK,PS,SRP,SRP2} for more information.

Parking function $(a_1,\cdots,a_n)$ can be redefined that its
increasing rearrangement $(b_1,\cdots,b_n)$ satisfies $b_i\leq i$.
 Pitman and  Stanley generalized the notion of parking functions
in \cite{PS}. Let ${\bf x}=(x_1,\cdots,x_n)$ be a sequence of
positive integers. The sequence $\alpha=(a_1,\cdots,a_n)$  is called
an ${\bf x}$-parking function if the non-decreasing rearrangement
$(b_1,\cdots,b_n)$ of $\alpha$ satisfies $b_i\leq x_1+\cdots +x_i$
for any $1\leq i\leq n$. Thus, the ordinary parking function is the
case ${\bf x}=(1,\cdots,1)$. By the determinant formula of
Gon\v{c}arove polynomials, Kung and Yan \cite{KY} obtained the
number of ${\bf x}$-parking functions for an arbitrary ${\bf x}$.
See also \cite{Y1,Y2,Y3} for the explicit formulas and properties
for some specified cases of ${\bf x}$.

An ${\bf x}$-parking function $(a_1,\cdots,a_n)$ is said to be
$k$-leading if $a_1=k$. Let $p_{n,k}$ denote the number of
$k$-leading ordinary parking functions of length $n$. Foata and
Riordan \cite{FR} derived the generating function for $p_{n,k}$
algebraically. Recently, Sen-peng Eu, Tung-shan Fu and Chun-Ju Lai
\cite{EFL} gave a combinatorial approach to the enumeration of
$(a,b,\cdots,b)$-parking functions by their leading terms.

Riordan \cite{R} told us the relations between ordered parking
functions and Dyck paths. Sen-peng Eu et al. \cite{EFY,ELY}
considered the problem of the enumerations of lattice paths with
flaws. It is natural to consider the problem of the enumerations of
preference sets with flaws. There are a few research results about
(ordered) $k$-flaw preference sets. There is a interesting facts.
Salmon and Cayley in 1849 established the classical configuration of
$27$ lines in a general cubic surface. Given a line $l$, the number
of lines, which are disjoint from, intersect or are equal to  $l$,
are $16$,$10$ and $1$, respectively, see \cite{henderson1911} for
the detail information. These number exactly are the number of $0$-,
$1$- and $2$-flaw preference sets of length $3$. Peter J Cameron et
al. \cite{Cam} indicate that there are some relations between
$k$-flaw preference sets and the context of hashing since data would
be lost. Also they counted the number of $k$-flaw preference sets
and calculate the asymptotic. In another joint work \cite{Huang},
using the methods different with \cite{Cam}, Po-Yi Huang et al.
studied $k$-flaw preference sets. In this paper, we focus on ordered
$k$-flaw preference sets. Let $\mathcal{OP}_{n,\geq k}$ denote the
set of ordered preference sets of length $n$ with at least $k$ flaws
and $op_{n,\geq k}=|\mathcal{OP}_{n,\geq k}|$. It is well known that
the number of solutions of the equation $ x_1+x_2+\cdots +x_k=n $ in
nonnegative integers is $n+k-1\choose{k-1}$. Hence, the number of
solutions of the equation $ x_1+x_2+\cdots +x_{n-k}=n+k $ in
nonnegative integers is $2n-1\choose{n-k-1}$. Let
$\mathcal{S}_{n,k}=\{(x_1,\cdots,x_{n-k})\mid x_1+x_2+\cdots
+x_{n-k}=n+k, x_i\in\mathbb{N}\}$. One of our main results is that
we obtain a bijection $\zeta$ between the sets $\mathcal{OP}_{n,\geq
k}$ and $\mathcal{S}_{n,k}$. With the benefit of the bijection , we
conclude that $op_{n,\geq k}={2n-1\choose{n-1-k}}$. Let
$\mathcal{OP}_{n,k}$($\mathcal{OP}_{n,\leq k}$) denote the set of
ordered preference sets with exact (at most) $k$ flaws and
$op_{n,k}=|\mathcal{OP}_{n, k}|$ ($op_{n,\leq
k}=|\mathcal{OP}_{n,\leq k}|$), respectively. By a simple
computation, $op_{n,k}=op_{n,\geq k}-op_{n,\geq k+1}$, we derive
that the number of order $k$-flaw preference sets is
$\frac{k+1}{n}{2n\choose{n-k-1}}$. Similarly, we also obtain
$op_{n,\leq k}={2n-1\choose{n-1}}-{2n-1\choose{n-k-2}}$.

 An $(n,\leq k)$-\emph{flaw path}
is a lattice path starting at $(0,0)$ and ending at $(2n,0)$ with
only two kinds of steps---\emph{rise step: $U=(1,1)$ }and \emph{fall
step: $D=(1,-1)$} lying on the line $y = -k$. Let
$\mathcal{D}_{n,\leq k}$ denote the set of $(n,\leq k)$-flaw paths.
If an $(n,\leq k)$-flaw path $P$ touches the line $y=-k$, then $P$
is called $(n, k)$-\emph{flaw path}. Let $\mathcal{D}_{n,k}$ denote
the set of $(n, k)$-flaw paths. We also find a bijection $\omega$
from the sets $\mathcal{OP}_{n,k}$ to $\mathcal{D}_{n,k}$, and this
gives the identity $op_{n,k}=\frac{k+1}{n}{2n\choose{n-k-1}}$ a
bijection proof. In fact, the mapping $\omega$ can be viewed as a
bijection from the sets $\mathcal{OP}_{n,\leq k}$ to
$\mathcal{D}_{n,\leq k}$.

Motivated by the work of Foata and Riordan in \cite{FR} as well as
Sen-Peng Eu et al. in \cite{EFL}, we investigate the problem of the
enumeration of ordered $k$-flaw preference sets with leading term
$m$. Let $\mathcal{OP}_{n,k}^m$ denote the set of ordered $k$-flaw
preference sets of length $n$ with leading term $m$ and
$op_{n,k}^m=|\mathcal{OP}_{n,k}^m|$. We find that
$op_{n,k}^{k+1}=\frac{k+2}{n+1}{2n-k-1\choose{n-k-1}}$ and
$op_{n,k}^m=\frac{2k-m+4}{2n-m}{2n-m\choose{n-k-2}}$ for any $m\leq
k$. Furthermore, let $\mathcal{OP}_{n,k,=l}^m$  denote the set of
ordered $k$-flaw preference sets of length $n$ with leading term $m$
and $\max\{a_i\mid i\in [n]\}=l$. Let
$op_{n,k,=l}^m=|\mathcal{OP}_{n,k,=l}^m|$. We conclude that
$op_{n,k,=l}^{k+1}=\frac{n-l+k+1}{n+l-k-1}{n+l-k-1\choose{l-k-1}}$
and
$op_{n,k,=l}^{m}=\frac{(n-l+2k-m+4)(n-m+k+1)-(l-k-2)}{(n-m+k+2)(n+l-m-1)}{n+l-m-1\choose{l-k-2}}$
for any $m\leq k$.

Also we are interested in the recurrence relations for the sequence
$op_{n,k}^m$. We find that for any $k\geq 0$, the sequence
$op_{n,k}^{k+1}$ satisfies the recurrence relation:
$op_{n+1,k+1}^{k+2}=op_{n,k}^{k+1}+op_{n,k+1}^{k+1}+op_{n,k+1}^{k+2}$;
for any  $1\leq m\leq k$, the sequence $op_{n,k}^m$ satisfies the
recurrence relation:
$op_{n,k}^m=\sum\limits_{j=m+1}^{k+1}\sum\limits_{i=1}^{n-k}c_iop_{n-i,k}^{j}$;
for any $k\geq 0$, the sequence $op_{n,k}$ satisfies the recurrence
relation:
$op_{n,k}=\sum\limits_{i=0}^{n-k-1}c_i[op_{n-i,k}-op_{n-i,k}^1]. $
We use $x$,$y$ and $z$ to mark the number of flaws, length and
leading term of a preference set, respectively. By these recurrence
relations, we easily obtain the generating functions for the
corresponding sequences. For any $k\geq 0$, the generating function
$\varphi_{k}(y)$ for $op_{n,k}^{k+1}$ satisfies the recurrence
relation: $ \varphi_{k+1}(y)=yC(y)\varphi_{k}(y)\text{ and
}\varphi_{0}(y)=y[C(y)]^2; $ for any $k\geq 0$ and $k\geq m\geq 1$,
the generating function $\rho_{m,k}(y)$ for $op_{n,k}^m$ satisfies
the recurrence relation: $\rho_{m,k}(y)=C(y)\rho_{m+1,k}(y)$ and $
\rho_{k,k}(y)=y^{k+2}[C(y)]^{k+4};$ for any $k\geq 0$, and the
generating function  $\phi_{k}(y)$ for  $op_{n,k}$ is
$\phi_k(y)=y^{k+1}[C(y)]^{2(k+1)}$. So, we easily compute the mean
and the variance of flaws of all $k$-flaw preference sets in
$\mathcal{OP}_{n,k}^m$ and $\mathcal{OP}_{n,k}$, respectively.

Recently, Postnikov and Shapiro \cite{postnikov2004} gave a new
generalization, building on work of Cori, Rossin and Salvy
\cite{cori2002}, the $G$-parking functions of a graph. For the
complete graph $G=K_{n+1}$, the defined functions in
\cite{postnikov2004} are exactly the classical parking functions.
So, in the future work, we will consider $k$-flaw $G$-parking
function.

We organize this paper as follows. The bijection $\zeta$ is given in
Section $2$. In Section $3$, we enumerate $k$-flaw preference sets
by the leading terms. In Section $4$, we obtain some the recurrence
relations of the sequence formed by some $k$-flaw preference sets.
In Section $5$, we derive the generating functions of the sequence
formed by $k$-flaw preference sets. In Section $6$, we compute the
mean and the variance of flaws of all $k$-flaw preference sets in
$\mathcal{OP}_{n,k}^m$ and $\mathcal{OP}_{n,k}$. We list the number
of $k$-flaw preference sets of length $n$ with leading term $m$
obtained by computer search for $n\leq 6$ in the Appendix.

\section{The bijection from $\mathcal{OP}_{n,\geq k}$ to $\mathcal{S}_{n,k}$}

 In this section, we will give the bijection $\zeta$ from the sets
$\mathcal{OP}_{n,\geq k}$ to $\mathcal{S}_{n,k}$.

Suppose $\alpha$ is an ordered $\hat{k}$-flaw preference set of
length $n$, Let ${\bf r}_\alpha=(r_1,\cdots,r_n)$ be the {\it
specification} of $\alpha$, i.e., $r_i=|\{j\mid a_j=i\}|.$ When
$\hat{k}=0$, Riordan \cite{R} proved the following lemma.
\begin{lemma}{\rm \cite{R}}\label{parkingfunctionr} The number of
parking functions of length $n$ is equals the $n$-th Catalan number
$c_n$. Furthermore, let $\alpha=(a_1,\cdots,a_n)$ be a parking
function of length $n$ and ${\bf r}_\alpha=(r_1,\cdots,r_n)$ the
{\it specification} of $\alpha$. Then for any $m\in[n]$,
$\sum\limits_{i=1}^{m}r_i-m\geq 0$ and
$\sum\limits_{i=1}^{n}r_i-n=0$.
\end{lemma}

Now, we consider the case with
$\alpha=(a_1,\cdots,a_n)\in\mathcal{OP}_{n,\geq k}$. Suppose
$\alpha$ is a ordered $\hat{k}$-flaw preference set, then
$\hat{k}\geq k$. Furthermore, we may suppose $\hat{k}\neq 0$ and
parking spaces $m_1,\cdots,m_{\hat{k}}$ are empty, then $r_{m_i}=0$
for $i\in[\hat{k}]$. Set $m_0=0$ and $m_{\hat{k}+1}=n+\hat{k}$. For
each $i\in[\hat{k}+1]$, let $T_{i}=\{j\mid m_{i-1}<a_j<m_{i}\}$ and
$u_i=|T_i|=m_i-m_{i-1}-1$.
\begin{lemma}\label{parkingfunction}
(1) Let $i\in[\hat{k}]$. For any $m\in[u_i]$ we have
$\sum\limits_{j=1}^{m}r_{m_{i-1}+j}-j\geq 0$ and
$\sum\limits_{j=1}^{u_i}r_{m_{i-1}+j}-u_i= 0$.

(2) For any $m\in[n-m_{\hat{k}}]$ we have
$\sum\limits_{j=1}^{m}r_{m_{\hat{k}}+j}-j\geq 0$ and
$\sum\limits_{j=1}^{n-m_{\hat{k}}}r_{m_{\hat{k}}+j}-(n-m_{\hat{k}})>
0$.
\end{lemma}
{\bf Proof.} (1) Let $\alpha_{T_i}=(a_{j_1},\cdots,a_{j_{u_i}})$ be
a subsequence of $\alpha$ determined by the subscripts in $T_i$ and
$\beta_i=(a_{j_1}-m_{i-1},\cdots,a_{j_{u_i}}-m_{i-1})$. Obviously,
$\beta_i$ is an ordered parking function of length $u_i$. By Lemma
\ref{parkingfunctionr}, we obtain the results as desired.

(2) Similar to (1), by Lemma \ref{parkingfunctionr}, for any
$m\in[n-m_{\hat{k}}]$ we have
$\sum\limits_{j=1}^{m}r_{m_{\hat{k}}+j}-j\geq 0$. Note $\hat{k}\neq
0$ and $u_{\hat{k}+1}=n+\hat{k}-m_{\hat{k}}>n-m_{\hat{k}}$. We have
$\sum\limits_{j=1}^{n-m_{\hat{k}}}r_{m_{\hat{k}}+j}-(n-m_{\hat{k}})>
0$.\hfill$\blacksquare$\\

Furthermore, Let $R=\{m_0,m_1,\cdots,m_{\hat{k}}\}$, $R+1=\{m+1\mid
m\in R\}$, $H=\{m_{\hat{k}-k+1},\cdots,m_{\hat{k}}\}$,
$H+1=\{m+1\mid m\in H\}$ and  $T=[n]\setminus H$. Obviously
$|T|=n-k$. Suppose $T=\{i_1,\cdots,i_{n-k}\}$ satisfying
$i_j<i_{j+1}$ for any $j\in[n-k-1]$ and ${\bf
r}_{\alpha}(T)=(r_{i_1},\cdots,r_{i_{n-k}})$ be a subsequence of
${\bf r}_{\alpha}$ determined by the subscripts in $T$.

For any $i_j\in T$, if $i_j\notin H+1$, then let $w_j=0$; if $i_j\in
H+1$, then let $w_j=|[i_{j-1},i_j]\cap H|$ when $j>1$ and
$w_1=|[1,i_1]\cap H|$. At last, let $x_{j}=r_{i_j}+w_j$ for all
$j=1,2,\cdots,n-k$ and $x=(x_1,\cdots,x_{n-k})$.
\begin{lemma}
For any $j\in[n-k]$, suppose $i_j\in T$. Let
$y(j,m)=\sum\limits_{l=0}^{m-1}x_{j+l}-m$ for all $m\in[n-k+1-j]$.
Then

(1) If $i_j\in (R+1)\setminus (H+1)$, then there exists a $\hat{m}$
satisfy $y(j,m)\geq 0$ for all $m\in[\hat{m}-1]$,  $y(j,\hat{m}-1)=
0$ and $y(j,\hat{m})<0$.

(2) If $i_j=m_s+1$ for some $s\in[\hat{k}-k+1,\hat{k}-1]$, then
$w_j=min\{y(j,m)\mid m\in[n-k+1-j]\}$ and $m_{s+1}-m_s-1=max\{m\mid
y(j,m)=w_j\}$

(3) If $i_j=\hat{k}+1$, then $w_j\leq min\{y(j,m)\mid
m\in[n-k+1-j]\}$.

(4) $\sum\limits_{i_j\in T\cap (H+1)}w_j=k$

(5) $\sum\limits_{i=1}^{n-k}x_i=n+k$.
\end{lemma}
{\bf Proof.} (1) By Lemma \ref{parkingfunction}, if $i_j=1$, then
$\hat{m}=m_1-1$; if $i_j=m_s+1$ for some $s\in[\hat{k}-k]$, then
$\hat{m}=m_{s+1}-m_{s}-1$.

(2) Since $i_j=m_s$ for some $s\in[\hat{k}-k+1,\hat{k}-1]$, by Lemma
\ref{parkingfunction}, we have $y(j,m)\geq w_j$ for all $m\in
[m_{s+1}-m_s-2]$ and $y(j,m_{s+1}-m_s-1)=w_j$. Note that
$x_{j+m_{s+1}-m_s-1}-1>0$. Hence, $w_j=min\{y(j,m)\mid
m\in[n-k+1-j]\}$ and $m_{s+1}-m_s-1=max\{m\mid y(j,m)=w_j\}$.

(3) Similar to (2), by Lemma \ref{parkingfunction}, we have $w_j\leq
min\{y(j,m)\mid m\in[n-k+1-j]\}$.

(4) When $i_j\in T\cap(H+1)$, since $w_j=|[i_{j-1},i_j]\cap H|$ if
$j>1$ and $w_1=|[1,i_1]\cap H|$, we have $\sum\limits_{i_j\in T\cap
(H+1)}w_j=k$.

(5) Clearly, $\sum\limits_{i_j\in T}r_{i_j}=n$. By (4), we have
$$\sum\limits_{i=1}^{n-k}x_i=\sum\limits_{i_j\in
T}(r_{i_j}+w_{i_j})=\sum\limits_{i_j\in
T}r_{i_j}+\sum\limits_{i_j\in
T\cap(H+1)}w_{i_j}=n+k.$$\hfill$\blacksquare$

\begin{example}
Take $n=9$ and $k=2$. Consider
$\alpha=(1,1,3,3,5,9,9,9,9)\in\mathcal{OP}_{9,\geq 2}$. Then ${\bf
r}_{\alpha}=(2,0,2,0,1,0,0,0,4)$, $R=\{6,7,8\}$, $H=\{7,8\}$, So,
$R+1=\{7,8,9\}$, $H+1=\{8,9\}$ and $T=\{1,2,3,4,5,6,9\}$. Hence
${\bf r}_{\alpha}(T)=(2,0,2,0,1,0,4)$, $w=(0,0,0,0,0,0,2)$. Finally,
we have ${x}=(2,0,2,0,1,0,6)$.
\end{example}

The above process give a mapping $\zeta$ from the sets
$\mathcal{OP}_{n,\geq k}$ to $\mathcal{S}_{n,k}$. In fact, we may
prove the mapping $\zeta$ is a bijection from the sets
$\mathcal{OP}_{n,\geq k}$ to $\mathcal{S}_{n,k}$.
\begin{theorem}\label{bijection}
There is a bijection $\zeta$ from the sets $\mathcal{OP}_{n,\geq k}$
to $\mathcal{S}_{n,k}$.
\end{theorem}
{\bf Proof.} To find $\zeta^{-1}$, let $x=(x_1,\cdots,x_{n-k})\in
\mathcal{S}_{n,k}$. For any $i\in[n-k]$ and $j\in[n-i+1]$, let
$y(i,j)=\sum\limits_{m=0}^{j-1}x_{i+m}-j$. If there is a
$m\in[n-i+2]$ such that $y(i,m+1)<0$, suppose $y(i,\widehat{m}+1)$
is the first negative number in the sequence
$y(i,1),\cdots,y(i,n-i+1)$, then let $u(i)=\widehat{m}$ and
$v(i)=0$. Otherwise, we have $y(i,j)\geq 0$ for all $j\in[n-i+1]$,
then let $\widehat{m}=max\{j'\mid y_{i,j'}\leq y_{i,j}\text{ for all
}j\in[n-i+1]\}$ and $u(i)=\widehat{m}$ and $v(i)=y(i,\widehat{m})$.
At last, let ${\bf u}=(u(1),\cdots,u(n-k))$ and ${\bf
v}=(v(1),\cdots,v(n-k))$.

Now, let
$$\begin{array}{ll}
\widehat{u}(1)=u(1)&\widehat{v}(1)=v(1)\\
\widehat{u}(i+1)=u(u(i)+1)&\widehat{v}(i+1)=v(u(i)+1),\end{array}
$$
and $s(i)=\sum\limits_{j=1}^{i}\widehat{u}(i)$ and
$t(i)=\sum\limits_{j=1}^{i}\widehat{v}(i)$, then there exists a
$\widetilde{m}$ satisfying $t(\widetilde{m})\leq k$ and
$t(\widetilde{m}+1)>k$. So, for any $i\in[\widetilde{m}]$, let ${\bf
r}_{i}=(r_{i,1},\cdots,r_{i,\widehat{u}_i+\widehat{v}_i})$ such that
$$
r_{1,j}=\left\{\begin{array}{lll} 0&\text{if}&1\leq j\leq
\widehat{v}(1)\\
x_{j-\widehat{v}(1)}-\widehat{v}(1)&\text{if}&j=\widehat{v}(1)+1\\
x_{j-\widehat{v}(1)}&\text{if}&\widehat{v}(1)+2\leq j\leq
\widehat{u}(1)+\widehat{v}(1)\end{array}\right.$$ for $2\leq i\leq
\widetilde{m}$,
$$
r_{i,j}=\left\{\begin{array}{lll} 0&\text{if}&1\leq j\leq
\widehat{v}(i)\\
x_{s(i-1)+j-v(i)}-\widehat{v}(i)&\text{if}&j=v(i)+1\\
x_{s(i-1)+j-v(i)}&\text{if}&v(i)+2\leq j\leq
\widehat{u}(i)+\widehat{v}(i)\\\end{array}\right.$$
 and ${\bf
r}_{\widetilde{m}+1}=(r_{\widetilde{m}+1,1},\cdots,r_{\widetilde{m}+1,n-s(\widetilde{m})-t(\widetilde{m})})$
such that
$$
r_{\widetilde{m}+1,j}=\left\{\begin{array}{lll} 0&\text{if}&1\leq
j\leq
k-t(\widetilde{m})\\
x_{s(\widetilde{m})+j-k+t(\widetilde{m})}-k+t(\widetilde{m})&\text{if}&j=k-t(\widetilde{m})+1\\
x_{s(\widetilde{m})+j-k+t(\widetilde{m})}&\text{if}&k-t(\widetilde{m})+2\leq
j\leq n-s(\widetilde{m})-t(\widetilde{m})\end{array}\right.$$

Finally, we obtain a specification ${\bf r}=({\bf r}_1,\cdots,{\bf
r}_{\widetilde{m}+1})$ of a preference set of length $n$ with at
least $k$ flaw. \hfill$\blacksquare$
\begin{example}
Let us consider the equation $x_1+\cdots+x_7=11$ and a solution
$x=(2,0,2,0,1,0,6)$. It is easy to obtain that $n=9$ and $k=2$. By
simple computations, we may get ${\bf u}=(6,0,3,0,1,0,5)$ and ${\bf
v}=(0,0,0,0,0,0,5)$. Clearly, $\widetilde{m}=1$. Hence, ${\bf
r}_1=(2,0,2,0,1,0)$ and ${\bf r}_2=(0,0,4)$. So, we have ${\bf
r}=(2,0,2,0,1,0,0,0,4)$. Finally, we construct a ordered preference
set $(1,1,3,3,5,9,9,9,9)$ of length $9$ with at least $2$ flaws.
\end{example}
With the benefit of the bijection in Theorem \ref{bijection}, we
obtain the following corollary.
\begin{corollary}\label{opatleastk-flaw}
$op_{n,\geq k}={2n-1\choose{n-1-k}}.$
\end{corollary}
\begin{example}Take $n=6$ and $k=3$. In Appendix, we can find $op_{6,\geq
3}=\sum\limits_{i=3}^{5}op_{6,i}=55={11\choose{2}}$.
\end{example}
\begin{remark} Note that $op_{n,k}=op_{n,\geq k}-op_{n,\geq k+1}$. By
Corollary \ref{opatleastk-flaw}, we have
$op_{n,k}=\frac{k+1}{n}{2n\choose{n-k-1}}.$ Let
$\mathcal{OP}_{n,\leq k}$ denote the set of preference sets of
length $n$ with at most $k$ flaws and $op_{n,\leq
k}=|\mathcal{OP}_{n,\leq k}|$. Obviously, $op_{n,\leq k}=op_{n,\geq
0}-op_{n,\geq k+1}$. Hence we have $op_{n,\leq
k}={2n-1\choose{n-1}}-{2n-1\choose{n-k-2}}$. In the rest of this
section, we will give the bijection proofs of the above results .
\end{remark}

Clearly, when $k=0$, an $(n,0)$-flaw path is an {\it n-Dyck path} as
well. We can also consider an $(n,k)$-flaw path $P$ a word of $2n$
letters using only $U$ and $D$. If a joint node is formed by a rise
step followed by a fall step, then this node is called a peak.

Let $\alpha=(a_1,\cdots,a_n)\in\mathcal{OP}_{n,k}$ and  ${\bf
r}_{\alpha}=(r_1,\cdots,r_n)$ be the specification  of $\alpha$. For
any $i\in[n]$, let $P_i=\underbrace{UU\cdots U}_{r_i}D$. So, we
obtain a lattice path $P=P_1\cdots P_n$ of length $2n$. Obviously,
$P$ is an $(n,k)$-flaw path and $P\in \mathcal{D}_{n,k}$. So, we
obtain a mapping from $\mathcal{OP}_{n,k}$ to $\mathcal{D}_{n,k}$,
and we  denote this mapping as $\omega$.
\begin{example}
Let $\alpha=(1,1,3,3,5,9,9,9,9)\in\mathcal{OP}_{9,3}$. Then ${\bf
r}_\alpha=(2,0,2,0,1,0,0,0,4)$. So,we may obtain a $(9,3)$-flaw path
$P=UUDDUUDDUDDDDUUUUD$.
\end{example}

In fact, we may prove the following theorem.
\begin{theorem}\label{bijectionlatticepath} There is a bijection $\omega$ from $\mathcal{OP}_{n,k}$ to
$\mathcal{D}_{n,k}$\end{theorem} {\bf Proof.} To find $\omega^{-1}$,
for any $P\in\mathcal{D}_{n,k}$, we view $P$ as a word of length
$2n$ using only $U$ and $D$. Our goal is to find a vector ${\bf r}$
such that ${\bf r}$ is the specification of a $k$-flaw preference
set of length $n$. By adding a symbol $*$ and replacing every peak
$UD$ with $U*$, we may obtain a new word $\widetilde{P}$. Now, we
suppose $\widetilde{P}=w_1w_2\cdots w_{2n}$. Let $T=\{j\mid
w_j=D\text{ or }w_j=*\}$. Clearly, $|T|=n$. Hence, we may suppose
$T=\{i_1,\cdots,i_n\}$ satisfying $i_s<i_{s+1}$. For any $j\in [n]$,
if $w_{i_j}=D$, then let $r_j=0$; if $w_{i_j}=*$, then $r_j=|\{s\mid
w_s=u\text{ and } i_{j-1}<s<i_j\}|$. So, we construct a vector ${\bf
r}=(r_1,\cdots,r_n)$. Since $\sum\limits_{i=1}^{n}r_i=n$, there is a
vector $\alpha=(a_1,\cdots,a_n)$ satisfying ${\bf r}_\alpha={\bf
r}$. Since $P$ touch the line $y=-k$, $\alpha$ is a $k$-flaw
preference set of length $n$.\hfill$\blacksquare$
\begin{example}
Let $P=UUDDUUDDUDDDDUUUUD$, clearly, $P\in\mathcal{D}_{9,3}$.
Furthermore,  we have $\widetilde{P}=UU*DUU*DU*DDDUUUU*$ and
$T=\{3,4,7,8,10,11,12,13,18\}$. Note that $i_1=3$ and
$w_{i_1}=w_3=*$, hence $r_1=2$. Since $i_2=4$ and $w_{i_2}=w_4=D$,
we have $r_2=0$. Finally, we may obtain a vector
$r=(2,0,2,0,1,0,0,0,4)$ and
$\alpha=(1,1,3,3,5,9,9,9,9)\in\mathcal{OP}_{9,3}$.
\end{example}

The mapping $\omega$ can be viewed as a bijection from the sets
$\mathcal{OP}_{n,\leq k}$ to $\mathcal{D}_{n,\leq k}$. So, we have
the following corollary.
\begin{corollary}
There is a bijection from the sets $\mathcal{OP}_{n,\leq k}$ to
$\mathcal{D}_{n,\leq k}$.
\end{corollary}
Let $d_{n,k}=|\mathcal{D}_{n,k}|$. It is well known that
$d_{n,k}=\frac{k+1}{n}{2n\choose{n-k-1}}$.

By Theorem \ref{bijectionlatticepath}, we immediately obtain the
following theorem.
\begin{theorem}\label{opk-flaw}For any $0\leq k\leq n-1$,
$op_{n,k}=\frac{k+1}{n}{2n\choose{n-k-1}}.$
\end{theorem}
\begin{example}
Take $n=6$ and $k=2$. In appendix, we find
$op_{6,2}=110=\frac{1}{2}{12\choose{3}}$.
\end{example}

\section{Enumerating ordered $k$-flaw preference sets by leading terms}
In this section, we will enumerate ordered $k$-flaw preference sets
by leading terms. First, let $\mathcal{OP}_{n,\geq k,\leq l}$ denote
the set of ordered preference sets $\alpha=(a_1,\cdots,a_n)$ of
length $n$ which have at least $k$ flaws and satisfy $a_i\leq l$ for
any $i\in[n]$. Let $op_{n,\geq k,\leq l}=|\mathcal{OP}_{n,\geq
k,\leq l}|$. It is easy to obtain that $l\geq n$ and $op_{n,\geq
k,\leq l}=0$ if $l\leq k$.
\begin{lemma}\label{atleastkflawmaxl} For any $n\geq l\geq k+1$,
$
op_{n,\geq k,\leq l}={n+l-1\choose{l-k-1}}.
$
\end{lemma}
{\bf Proof.} For any $n\geq l\geq k+1$, be Theorem \ref{bijection},
$op_{n,\geq k,\leq l}$ is equals the number of the solutions of
equation $ x_1+x_2+\cdots+x_{l-k}=n+k $ in nonnegative integers.
Hence, $ op_{n,\geq k,\leq l}={n+l-1\choose{l-k-1}}.
$ \hfill$\blacksquare$\\

Now, for any $m\leq l$, let $\mathcal{OP}_{n,\geq k,\leq l}^m$
denote the set of ordered preference sets
$\alpha=(a_1,\cdots,a_n)\in\mathcal{OP}_{n,\geq k,\leq l}$ with
leading term $m$. Let $op_{n,\geq k,\leq l}^m=|\mathcal{OP}_{n,\geq
k,\leq l}^m|$.

\begin{lemma}\label{atleastkflawmaxlleaingk+1}For any $n\geq l\geq
k+1$, $op_{n,\geq k,\leq l}^{k+1}={n+l-k-2\choose{l-k-1}}.$
\end{lemma}
{\bf Proof.} Suppose $\alpha=(a_1,\cdots,a_n)\in\mathcal{OP}_{n,\geq
k,\leq l}^{k+1}$, let $\beta=(a_1-k,\cdots,a_n-k)$, then
$\beta\in\mathcal{OP}_{n,\geq 0,\leq l-k}^1$. Conversely, for any
$\beta=(b_1,\cdots,b_n)\in\mathcal{OP}_{n,\geq 0,\leq l-k}^1$, let
$\alpha=(b_1+k,\cdots,b_n+k)$, then $\alpha\in\mathcal{OP}_{n,\geq
k,\leq l}^{k+1}$. So, $op_{n,\geq k,\leq l}^{k+1}=op_{n,\geq 0,\leq
l}^1$

By Theorem \ref{bijection}, $op_{n,\geq 0,\leq l}^1$ is equals the
number of the solutions of equation
$$\left\{\begin{array}{l} x_1+x_2+\cdots+x_{l}=n\\
x_1\geq 1,x_2\geq 0,\cdots,x_l\geq 0,
\end{array}\right.$$
in nonnegative integers. Let $y_1=x_1-1$ and $y_i=x_i$ for all
$i\geq 2$. So, $op_{n,\geq 0,\leq l}^1$ is equals the number of the
solutions of equation $ y_1+y_2+\cdots+y_{l}=n-1  $ in nonnegative
integers. Hence, we have $op_{n,\geq k,\leq
l}^{k+1}={n+l-k-2\choose{l-k-1}}$ for any $n\geq l\geq k+1$.
\hfill$\blacksquare$

\begin{lemma}\label{allestkflawmaxlleadingm}For any $n\geq l\geq k+1$ and $k\geq m$,
$op_{n,\geq k,\leq l}^m={n+l-m-1\choose{l-k-2}}.$
\end{lemma}
{\bf Proof.} Suppose $\alpha=(a_1,\cdots,a_n)\in\mathcal{OP}_{n,\geq
k,\leq l}^m$, let $\beta=(a_1-m+1,\cdots,a_n-m+1)$, then
$\beta\in\mathcal{OP}_{n,\geq k-m+1,\leq l-m+1}^1$. Conversely, for
any $\beta=(b_1,\cdots,b_n)\in\mathcal{OP}_{n,\geq k-m+1,\leq
l-m+1}^1$, let $\alpha=(b_1+m-1,\cdots,b_n+m-1)$, then
$\alpha\in\mathcal{OP}_{n,\geq k,\leq l}^m$. Hence, $op_{n,\geq
k,\leq l}^m=op_{n,\geq k-m+1,\leq l-m+1}^1$

For any
$\alpha=(a_1,\cdots,a_n)\in\bigcup\limits_{i=2}^{l-m+1}\mathcal{OP}_{n,\geq
k-m+1,\leq l-m+1}^i$, let $\beta=(a_1-1,\cdots,a_n-1)$, then
$\beta\in\mathcal{OP}_{n,\geq k-m,l-m}$. Conversely, for any
$\beta=(b_1,\cdots,b_n)\in\mathcal{OP}_{n,\geq k-m,l-m}$, let
$\alpha=(b_1+1,\cdots,b_n+1)$, then
$\alpha\in\bigcup\limits_{i=2}^{l-m+1}\mathcal{OP}_{n,\geq
k-m+1,\leq l-m+1}^i$. So, $\sum\limits_{i=2}^{l-m+1}op_{n,\geq
k-m+1,\leq l-m+1}^i=op_{n,\geq k-m,\leq l-m}$.

Hence, by Lemma \ref{atleastkflawmaxl}, we have
\begin{eqnarray*} op_{n,\geq
k,\leq l}^m&=&op_{n,\geq k-m+1,\leq l-m+1}^1\\
&=&op_{n,\geq
k-m+1,\leq l-m+1}-\sum\limits_{i=2}^{l-m+1}op_{n,\geq
k-m+1,\leq l-m+1}^m\\
&=&op_{n,\geq k-m+1,\leq l-m+1}-op_{n,\geq k-m,\leq l-m}\\
&=&{n+l-m\choose{l-k-1}}-{n+l-m-1\choose{l-k-1}}\\
&=&{n+l-m-1\choose{l-k-2}}.
\end{eqnarray*}\hfill$\blacksquare$\\

Let $\mathcal{OP}_{n,k,\leq l}^m$ denote the set of ordered $k$-flaw
preference sets $\alpha=(a_1,\cdots,a_n)$ of length $n$ which
satisfy $a_1=m$ and $a_i\leq l$ for any $i\in[n]$. Let $op_{n,
k,\leq l}^m=|\mathcal{OP}_{n, k,\leq l}^m|$. By Lemmas
\ref{atleastkflawmaxlleaingk+1} and \ref{allestkflawmaxlleadingm},
we obtain the following two corollaries.

\begin{corollary}\label{kflawmaxlleadk+1}For any $n\geq l\geq k+1$,
$op_{n,k,\leq l}^{k+1}=\frac{n-l+k+2}{n+1}{n+l-k-1\choose{l-k-1}}.$
\end{corollary}
{\bf Proof.} Note that $op_{n,k,\leq l}^{k+1}=op_{n,\geq k,\leq
l}^{k+1}-op_{n,\geq k+1,\leq l}^{k+1}$. By Lemma
\ref{atleastkflawmaxlleaingk+1}, we have
\begin{eqnarray*}op_{n,k,\leq l}^{k+1}&=&{n+l-k-2\choose{l-k-1}}-{n+l-k-2\choose{l-k-3}}\\
&=&\frac{n-l+k+2}{n+1}{n+l-k-1\choose{l-k-1}}.\end{eqnarray*}\hfill$\blacksquare$

\begin{corollary}\label{kflawmaxlleadm}For any $k\geq m\geq 1$ and
$n\geq l\geq k+1$, $op_{n,k,\leq
l}^{m}=\frac{n-l+2k-m+4}{n-m+k+2}{n+l-m-1\choose{l-k-2}}.$
\end{corollary}
{\bf Proof.} Note that $op_{n,k,\leq l}^{m}=op_{n,\geq k,\leq
l}^{m}-op_{n,\geq k+1,\leq l}^{m}$. By Lemma
\ref{allestkflawmaxlleadingm}, we have
\begin{eqnarray*}op_{n,k,\leq l}^{m}&=&{n+l-m-1\choose{l-k-2}}-{n+l-m-1\choose{l-k-3}}\\
&=&\frac{n-l+2k-m+4}{n-m+k+2}{n+l-m-1\choose{l-k-2}}.\end{eqnarray*}\hfill$\blacksquare$

For any $\alpha\in \mathcal{OP}_{n,k}^m$, it is easy to check that
parking spaces $1,2,\cdots,m-1$ are empty. Hence, $op_{n,k}^m=0$ if
$m\geq k+2$. Obviously, $op_{n,k}=0$ if $n\leq k$. So, we always
suppose $n\geq k+1\geq m$.

First, we consider the case with $m=k+1$.

\begin{theorem}\label{m-1flawleadingm}
Let $n\geq k+1$. Then
$op_{n,k}^{k+1}=\frac{k+2}{n+1}{2n-k-1\choose{n-k-1}}.$
\end{theorem}
{\bf Proof.} Take $l=n$ in Corollary
\ref{kflawmaxlleadk+1}.\hfill$\blacksquare$
\begin{example}Take $n=6$ and $k=2$. In Appendix, we find
$op_{6,2}^3=48=\frac{4}{7}{9\choose{3}}.$
\end{example}
\begin{theorem}\label{kflawleadingm}
Let $k\geq m\geq 1$ and $n\geq k+1$. Then
$op_{n,k}^m=\frac{2k-m+4}{2n-m}{2n-m\choose{n-k-2}}.$
\end{theorem}
{\bf Proof.} Take $l=n$ in Corollary
\ref{kflawmaxlleadm}.\hfill$\blacksquare$
\begin{example}Take $n=6$ , $k=2$ and $m=2$. In Appendix, we find
$op_{6,2}^2=27=\frac{3}{5}{10\choose{2}}.$
\end{example}
\begin{theorem}\label{kflawmaxlleadk+1existl}For any
$n\geq l\geq k+1$,
$op_{n,k,=l}^{k+1}=\frac{n-l+k+1}{n+l-k-1}{n+l-k-1\choose{l-k-1}}.$
\end{theorem}
{\bf Proof.} Note that $op_{n,k,=l}^{k+1}=op_{n, k,\leq
l}^{k+1}-op_{n,k,\leq l-1}^{k+1}$. By Corollary
\ref{kflawmaxlleadk+1}, we have
\begin{eqnarray*}op_{n,k,=l}^{k+1}&=&\frac{n-l+k+2}{n+1}{n+l-k-1\choose{l-k-1}}-\frac{n-l+k+3}{n+1}{n+l-k-2\choose{l-k-2}}\\
&=&\frac{n-l+k+1}{n+l-k-1}{n+l-k-1\choose{l-k-1}}.\end{eqnarray*}\hfill$\blacksquare$
\begin{theorem}\label{kflawmaxlleadmexistl}For any $k\geq m\geq 1$ and
$n\geq l\geq k+1$,
\begin{eqnarray*}op_{n,k,=l}^{m}=\frac{(n-l+2k-m+4)(n-m+k+1)-(l-k-2)}{(n-m+k+2)(n+l-m-1)}{n+l-m-1\choose{l-k-2}}.\end{eqnarray*}
\end{theorem}
{\bf Proof.} Note that $op_{n,k,=l}^{m}=op_{n,k,\leq
l}^{m}-op_{n,k,\leq l-1}^{m}$. By Corollary \ref{kflawmaxlleadm}, we
have
\begin{eqnarray*}op_{n,k,=l}^{m}&=&\frac{n-l+2k-m+4}{n-m+k+2}{n+l-m-1\choose{l-k-2}}-\frac{n-l+2k-m+5}{n-m+k+2}{n+l-m-2\choose{l-k-3}}\\
&=&\frac{(n-l+2k-m+4)(n-m+k+1)-(l-k-2)}{(n-m+k+2)(n+l-m-1)}{n+l-m-1\choose{l-k-2}}.\end{eqnarray*}\hfill$\blacksquare$

\section{The recurrence relation of ordered $k$-flaw preference sets}
In this section, we will give some recurrence relations of the
sequence formed by ordered $k$-flaw preference sets of length $n$
with leading term $m$.

First, we consider the case with $m=k+1$.
\begin{theorem}\label{kflawleadingk+1recurrence}For any $k\geq 0$,
the sequence formed by $k$-flaw preference sets of length $n$ with
leading term $k+1$ satisfies the following recurrence relation:
\begin{equation}op_{n+1,k+1}^{k+2}=op_{n,k}^{k+1}+op_{n,k+1}^{k+1}+op_{n,k+1}^{k+2}.\end{equation}
\end{theorem}
{\bf Proof.} Let $A$ denote the set of preference sets of length $n$
with at most $k+1$ flaws and leading term no less than $k+1$.
Obviously, for any $\beta=(b_1,\cdots,b_n)\in A$, we have $b_1=k+1$
or $k+2$, hence,
$A=\mathcal{OP}_{n,k}^{k+1}\cup\mathcal{OP}_{n,k+1}^{k+1}\cup\mathcal{OP}_{n,k+1}^{k+2}$
 and $|A|=op_{n,k}^{k+1}+op_{n,k+1}^{k+1}+op_{n,k+1}^{k+2}$.

Now, for any
$\alpha=(a_1,\cdots,a_{n+1})\in\mathcal{OP}_{n+1,k+1}^{k+2}$, let
$\beta=(a_2-1,\cdots,a_{n+1}-1)$, then $\beta\in A$; conversely, for
any $\beta=(b_1,\cdots,b_{n})\in A$, let
$\alpha=(k+2,b_1+1,\cdots,b_n+1)$, then
$\alpha\in\mathcal{OP}_{n+1,k+1}^{k+2}$. So,
$op_{n+1,k+1}^{k+2}=|\mathcal{OP}_{n+1,k+1}^{k+2}|=|A|=op_{n,k}^{k+1}+op_{n,k+1}^{k+1}+op_{n,k+1}^{k+2}.$
\hfill$\blacksquare$
\begin{example}Take $n=5$ and $k=1$. In Appendix, we can find $op_{6,2}^3=48$,
$op_{5,1}^2=28$,$op_{5,2}^2=6$ and $op_{5,2}^3=14$. Then
$op_{6,2}^3=op_{5,1}^2+op_{5,2}^2+op_{5,2}^3$.
\end{example}
\begin{remark}In fact, by Theorems \ref{m-1flawleadingm} and
\ref{kflawleadingm}, simple computation also implies the recurrence
relation in Theorem \ref{kflawleadingk+1recurrence}.
\end{remark}
\begin{theorem}\label{kflawleadingk+1max=lrecurrence}For any $k\geq 0$,
the sequence formed by $k$-flaw preference sets of length $n$ with
leading term $k+1$ and $\max\{a_i\mid i\in[n]\}=l$ satisfies the
following recurrence relation:
\begin{equation}op_{n+1,k+1,l+1}^{k+2}=op_{n,k,=l}^{k+1}+op_{n,k+1,=l}^{k+1}+op_{n,k+1,=l}^{k+2}.\end{equation}
\end{theorem}
{\bf Proof.} Let $\widetilde{A}$ denote the set of parking functions
of length $n$ with at most $k+1$ flaws and leading term no less than
$k+1$ and $\max\{a_i\mid i\in[n]\}=l$. Obviously, for any
$\beta=(b_1,\cdots,b_n)\in \widetilde{A}$, we have $b_1=k+1$ or
$k+2$, hence,
$\widetilde{A}=\mathcal{OP}_{n,k,=l}^{k+1}\cup\mathcal{OP}_{n,k+1,=l}^{k+1}\cup\mathcal{OP}_{n,k+1,=l}^{k+2}$
 and $|\widetilde{A}|=op_{n,k,=l}^{k+1}+op_{n,k+1,=l}^{k+1}+op_{n,k+1,=l}^{k+2}$.

Now, for any
$\alpha=(a_1,\cdots,a_{n+1})\in\mathcal{OP}_{n+1,k+1,l+1}^{k+2}$,
let $\beta=(a_2-1,\cdots,a_{n+1}-1)$, it is easy to check that
$\beta\in \widetilde{A}$; conversely, for any
$\beta=(b_1,\cdots,b_{n})\in \widetilde{A}$, let
$\alpha=(k+2,b_1+1,\cdots,b_n+1)$, then
$\alpha\in\mathcal{OP}_{n+1,k+1,l+1}^{k+2}$. So,
\begin{eqnarray*}op_{n+1,k+1,l+1}^{k+2}=|\mathcal{OP}_{n+1,k+1,l+1}^{k+2}|=|\widetilde{A}|=op_{n,k,=l}^{k+1}+op_{n,k+1,=l}^{k+1}+op_{n,k+1,=l}^{k+2}.\end{eqnarray*}
\hfill$\blacksquare$

 For the case with $1\leq m\leq k$, we obtain the following results.
 \begin{theorem}\label{kflawleadingmrecurrence} Let $1\leq m\leq k$.
The sequence formed by $k$-flaw preference sets of length $n$ with
leading term $m$ satisfies the following recurrence relation:
\begin{equation}op_{n,k}^m=\sum\limits_{j=m+1}^{k+1}\sum\limits_{i=1}^{n-k}c_iop_{n-i,k}^{j}.\end{equation}
 \end{theorem}
{\bf Proof.} Let $k\geq m\geq 1$,
$\alpha=(a_1,\cdots,a_n)\in\mathcal{OP}_{n,k}^m$ and
$\beta=(a_1-m+1,\cdots,a_n-m+1)$, then
$\beta\in\mathcal{OP}_{n,k-m+1,\leq n-m+1}^1$. Thus, we suppose the
$(i+1)$-th parking space is the first empty parking space, let
$S=\{j\mid a_j-m+1\leq i\}$ and $T=\{j\mid i+2\leq a_j-m+1\leq
n-m+1\}$, then $|S|=i$ and $|T|=n-i$. Let $\beta_{S}$ and
$\beta_{T}$ be the sequences of $\beta$ determined by the subscripts
in $S$ and $T$, respectively. We have $\beta_{S}\in\mathcal{OP}_i$
and $\beta_{T}\in\bigcup\limits_{j=
i+2}^{k-m+i+2}\mathcal{OP}_{n-i,k-m+i+1,\leq n-m+1}^j$.

Note that $1\leq i\leq n-k$. Let $c_n$ be $n$-th Catalan number .
Also it is easy to check $\bigcup\limits_{j=
i+2}^{k-m+i+2}\mathcal{OP}_{n-i,k-m+i+1,\leq
n-m+1}^j=\bigcup\limits_{j= m+1}^{k+1}\mathcal{OP}_{n-i,k}^j$.
Therefore, there are $c_i$ and
$\sum\limits_{j=m+1}^{k+1}op_{n-i,k}^{j}$ possibilities for
$\beta_{S}$ and $\beta_{T}$, respectively. Finally, we obtain the
following equation
\begin{eqnarray*}op_{n,k}^m=\sum\limits_{i=1}^{n-k}c_i\sum\limits_{j=m+1}^{k+1}op_{n-i,k}^{j}=\sum\limits_{j=m+1}^{k+1}\sum\limits_{i=1}^{n-k}c_iop_{n-i,k}^{j}.\end{eqnarray*}
\hfill$\blacksquare$
\begin{example}Take $n=6$, $k=2$ and $m=1$. By the data in Appendix,
it is easy to check
$op_{6,2}^1=\sum\limits_{j=2}^{3}\sum\limits_{i=1}^{4}c_iop_{6-i,2}^{j}$.
\end{example}
 \begin{theorem}\label{kflawleadingmmax=lrecurrence} Let $1\leq m\leq k$.
The sequence formed by $k$-flaw preference sets of length $n$ with
leading term $m$ and $\max\{a_i\mid i\in [n]\}=l$ satisfies the
following recurrence relation:
\begin{equation}op_{n,k,=l}^m=\sum\limits_{j=m+1}^{k+1}\sum\limits_{i=1}^{n-k}c_iop_{n-i,k,l-i}^{j}.\end{equation}
 \end{theorem}
{\bf Proof.} Let $k\geq m\geq 1$,
$\alpha=(a_1,\cdots,a_n)\in\mathcal{OP}_{n,k,=l}^m$ and
$\beta=(a_1-m+1,\cdots,a_n-m+1)$, then
$\beta\in\mathcal{OP}_{n,k-m+1,l-m+1}^1$. Thus, we suppose the
$(i+1)$-th parking space is the first empty parking space, let
$S=\{j\mid a_j-m+1\leq i\}$ and $T=\{j\mid i+2\leq a_j-m+1\leq
l-m+1\}$, then $|S|=i$ and $|T|=n-i$. Let $\beta_{S}$ and
$\beta_{T}$ be the sequences of $\beta$ determined by the subscripts
in $S$ and $T$, respectively. We have $\beta_{S}\in\mathcal{OP}_i$
and $\beta_{T}\in\bigcup\limits_{j=
i+2}^{k-m+i+2}\mathcal{OP}_{n-i,k-m+i+1,l-m+1}^j$.

Note that $1\leq i\leq n-k$. Let $c_n$ be $n$-th Catalan number .
Also it is easy to check $\bigcup\limits_{j=
i+2}^{k-m+i+2}\mathcal{OP}_{n-i,k-m+i+1,l-m+1}^j=\bigcup\limits_{j=
m+1}^{k+1}\mathcal{OP}_{n-i,k,l-i}^j$. Therefore, there are $c_i$
and $\sum\limits_{j=m+1}^{k+1}op_{n-i,k,l-i}^{j}$ possibilities for
$\beta_{S}$ and $\beta_{T}$, respectively. Finally, we obtain the
following equation
\begin{eqnarray*}op_{n,k,=l}^m=\sum\limits_{i=1}^{n-k}c_i\sum\limits_{j=m+1}^{k+1}op_{n-i,k,l-i}^{j}=\sum\limits_{j=m+1}^{k+1}\sum\limits_{i=1}^{n-k}c_iop_{n-i,k,l-i}^{j}.\end{eqnarray*}
\hfill$\blacksquare$

\begin{theorem}\label{kflawrecurrence}
Let $k\geq 0$. The sequence formed by $k$-flaw preference sets of
length $n$ satisfies the following recurrence relation:
\begin{eqnarray*}op_{n,k}=\sum\limits_{i=0}^{n-k-1}c_i[op_{n-i,k}-op_{n-i,k}^1].
\end{eqnarray*}
\end{theorem}
{\bf Proof.} Let $k\geq 0$,
$\alpha=(a_1,\cdots,a_n)\in\mathcal{OP}_{n,k}$. We suppose the
$(i+1)$-th parking space is the first empty parking space, let
$S=\{j\mid a_j\leq i\}$ and $T=\{j\mid i+2\leq a_j\leq n\}$, then
$|S|=i$ and $|T|=n-i$. Let $\alpha_{S}$ and $\alpha_{T}$ be the
sequences of $\alpha$ determined by the subscripts in $S$ and $T$,
respectively. Furthermore, suppose
$\alpha_{T}=(a_{i_1},\cdots,a_{i_{n-i}})$, let
$\beta_{T}=(a_{i_1}-i,\cdots,a_{i_{n-i}}-i)$, We have
$\alpha_{S}\in\mathcal{OP}_i$ and
$\beta_{T}\in\mathcal{OP}_{n-i,k}\setminus\mathcal{OP}_{n-i,k}^1$.

Note that $0\leq i\leq n-k-1$. Let $c_n$ be $n$-th Catalan number .
Therefore, there are $c_i$ and $op_{n-i,k}-op_{n-i,k}^1$
possibilities for $\alpha_{S}$ and $\beta_{T}$, respectively. We
obtain the following equation
$op_{n,k}=\sum\limits_{i=0}^{n-k-1}c_i[op_{n-i,k}-op_{n-i,k}^1].
$\hfill$\blacksquare$
\begin{example}Take $n=6$ and $k=2$. By the data in Appendix,
it is easy to check
$op_{6,2}=\sum\limits_{i=0}^{3}c_i[op_{6-i,2}-op_{6-i,2}^1]. $.
\end{example}
\section{Generating functions of ordered $k$-flaw preference sets}
In this section, we want to obtain the generating functions of the
sequences formed by ordered $k$-flaw preference sets of length $n$
with leading term $m$.

 Let $\psi(x,y)=\sum\limits_{n\geq
1}\sum\limits_{k=0}^{n-1}op_{n,\geq k}x^ky^n$ be the generating
function for $op_{n,\geq k}$.  To obtain $\psi(x,y)$, we need the
following lemma.
\begin{lemma}\label{B(y)} For any $n\geq 0$, let $b_{n}={2n\choose{n}}$ and
$B(x)=\sum\limits_{n\geq 0}b_nx^n$ be the generating function of the
sequence $b_1,b_2,\cdots$.Then $B(x)=\frac{1}{\sqrt{1-4x}}.$
\end{lemma}
{\bf Proof.}  Note that $b_n=(n+1)c_n$.  Hence,
 $
 B(x)=\frac{d(xC(x))}{dx}=\frac{1}{\sqrt{1-4x}}
 $\hfill$\blacksquare$
\begin{theorem}\label{opatleastk-flawgenerating}
Let $n\geq 1$. The generating function $\psi_{n}(x)$ for the
sequence formed by ordered parking functions length $n$ with at
least $k$ flaws satisfies the following recurrence relation:
\begin{eqnarray*}
2x\psi_n(x)=2(1+x)^2\psi_{n-1}(x)+(x-1){2n-2\choose{n-1}}.
\end{eqnarray*}

\end{theorem}
{\bf Proof.} By Lemma \ref{opatleastk-flaw}, using the equation
${n\choose m}={n-1\choose{m}}+{n-1\choose{m-1}}$, we obtain the
following recurrence relation.
\begin{eqnarray*}
\psi_n(x)=\frac{(1+x)^2}{x}\psi_{n-1}(x)+\frac{1}{2}(1-\frac{1}{x}){2n-2\choose{n-1}}.
\end{eqnarray*}\hfill$\blacksquare$
\begin{corollary}\label{atleastkflawgenerating}
Let  $\psi(x,y)$ be the generating function for the sequence formed
by ordered parking functions length $n$ with at least $k$ flaws.
Then
\begin{eqnarray*}\psi(x,y)
=\frac{y}{2(x-y(1+x)^2)}\left[\frac{x-1}{\sqrt{1-4y}}+x+1\right].\end{eqnarray*}
\end{corollary}
{\bf Proof.} By Theorem \ref{opatleastk-flawgenerating}, we have
\begin{eqnarray*} 2x\sum\limits_{n\geq
2}\psi_n(x)y^n&=&2(1+x)^2\sum\limits_{n\geq
2}\psi_{n-1}(x)y^n+(x-1)\sum\limits_{n\geq
2}{2n-2\choose{n-1}}y^n\\
&=&2y(1+x)^2\sum\limits_{n\geq
1}\psi_{n}(x)y^n+(x-1)y\sum\limits_{n\geq 1}{2n\choose{n}}y^n.
\end{eqnarray*}
Note $\psi_1(x)=1$. By Lemma \ref{B(y)}, we have
\begin{eqnarray*}
2x[\psi(x,y)-y]=2y(1+x)^2\psi(x,y)+(x-1)y\left[B(y)-1\right].
\end{eqnarray*}
Solving this equation, we have
\begin{eqnarray*}\psi(x,y)=\frac{y}{2(x-y(1+x)^2)}\left[\frac{x-1}{\sqrt{1-4y}}+x+1\right].\end{eqnarray*}\hfill$\blacksquare$\\

Let $\varphi(x,y)=\sum\limits_{n\geq
1}\sum\limits_{k=0}^{n-1}op_{n,k}^{k+1}x^{k}y^n=\sum\limits_{k\geq
0}\sum\limits_{n\geq k+1}op_{n,k}^{k+1}y^nx^{k}$ be the generating
function for $op_{n,k}^{k+1}$. Furthermore, let
$\varphi_{k}(y)=\sum\limits_{n\geq k+1}op_{n,k}^{k+1}y^n$ for any
$k\geq 0$, then $\varphi(x,y)=\sum\limits_{k\geq
0}\varphi_{k}(y)x^{k}$.
\begin{theorem}\label{generatingkflawleadingk+1}Let $k\geq 0$. The generating function $\varphi_{k}(y)$ of the sequence formed by
$k$-flaw preference sets length $n$ with leading term $k+1$
satisfies the following recurrence relation:
\begin{eqnarray*}
\varphi_{k+1}(y)=yC(y)\varphi_{k}(y)\text{ and
}\varphi_{0}(y)=y[C(y)]^2.
\end{eqnarray*}
\end{theorem}
{\bf Proof.} Obviously, when $k=0$, we have
$\varphi_{0}(y)=C(y)-1=y[C(y)]^2$ since $n\geq 1$.

For any $k\geq 1$, taking $m=k$ in Equations $(1)$ and $(3)$, we
have
\begin{eqnarray*}op_{n+1,k+1}^{k+2}=op_{n,k}^{k+1}+\sum\limits_{i=1}^{n-k-1}c_iop_{n-i,k+1}^{k+2}+op_{n,k+1}^{k+2}.\end{eqnarray*}
Hence,
\begin{eqnarray*}\sum\limits_{n\geq k+2}op_{n+1,k+1}^{k+2}y^n=\sum\limits_{n\geq k+2}op_{n,k}^{k+1}y^n+\sum\limits_{n\geq k+2}\sum\limits_{i=1}^{n-k-1}c_iop_{n-i,k+1}^{k+2}y^n+\sum\limits_{n\geq k+2}op_{n,k+1}^{k+2}y^n.\end{eqnarray*}
So,
\begin{eqnarray*}[\varphi_{k+1}(y)-op_{k+2,k+1}^{k+1}y^{k+2}]y^{-1}=\varphi_{k}(y)-op_{k+1,k}^{k}y^{k+1}
+[C(y)-1]\varphi_{k+1}(y)+\varphi_{k+1}(y).\end{eqnarray*} Note that
$op_{k+2,k+1}^{k+1}=op_{k+1,k}^{k}=1$ and $\frac{1}{C(y)}=1-yC(y)$.
Therefore,
$\varphi_{k+1}(y)=yC(y)\varphi_{k}(y)$.\hfill$\blacksquare$
\begin{corollary}Let $\varphi(x,y)$ be the generating function  of the sequence formed by
$k$-flaw preference sets length $n$ with leading term $k+1$. Then
\begin{eqnarray*}\varphi(x,y)
=\frac{y[C(y)]^2}{1-xyC(y)}.\end{eqnarray*}\end{corollary} {\bf
Proof.} By Theorem \ref{generatingkflawleadingk+1}, it is easy to
obtain that $\varphi_{k}(y)=y^{k+1}[C(y)]^{k+2}$.Note that
$\varphi(x,y)=\sum\limits_{k\geq 0}\varphi_{k}(y)x^{k}$. Hence,
\begin{eqnarray*}\varphi(x,y)=\sum\limits_{k\geq 0}\varphi_k(y)x^{k}
=\sum\limits_{k\geq 0}y^{k+1}[C(y)]^{k+2}x^{k}
=\frac{y[C(y)]^2}{1-xyC(y)}.  \end{eqnarray*}\hfill$\blacksquare$\\

Let $\rho(x,y,z)=\sum\limits_{m\geq 1}\sum\limits_{n\geq
m+1}\sum\limits_{k= m}^{n-1}op_{n,k}^mx^ky^nz^{m}=\sum\limits_{m\geq
1}\sum\limits_{k\geq m}\sum\limits_{n\geq k+1}op_{n,k}^my^nx^{k}z^m$
be the generating function of the sequence $op_{n,k}^m$ with $k\geq
m$. Let $\rho_m(x,y)=\sum\limits_{k\geq m}\sum\limits_{n\geq
k+1}op_{n,k}^my^nx^{k}$ and $\rho_{m,k}(y)=\sum\limits_{n\geq
k+1}op_{n,k}^my^n$ for any $k\geq m$, then
$\rho(x,y,z)=\sum\limits_{m\geq 1}\rho_m(x,y)z^m$ and
$\rho_m(x,y)=\sum\limits_{k\geq m}\rho_{m,k}(y)x^{k}$.
\begin{theorem}\label{generatingrecurrencemk}
Let $k\geq 0$ and $k\geq m\geq 1$. The generating function
$\rho_{m,k}(y)$ of the sequence formed by $k$-flaw preference sets
length $n$ with leading term $m$ satisfies the following recurrence
relation:
$$\left\{\begin{array}{l}\rho_{m,k}(y)=C(y)\rho_{m+1,k}(y)\\
\rho_{k,k}(y)=y^{k+2}[C(y)]^{k+4}.\end{array}\right.$$
\end{theorem}
{\bf Proof.}  By Equation $(3)$, we have $\sum\limits_{n\geq
k+1}op_{n,k}^my^n =\sum\limits_{j=m+1}^{k+1}\sum\limits_{n\geq
k+1}\sum\limits_{i=1}^{n-k}c_iop_{n-i,k}^{j}y^n.$

 {\bf Case 1.}
$k=m$

Note that $C(y)-1=y[C(y)]^2$. By Theorem
\ref{generatingkflawleadingk+1},  we have \begin{eqnarray*}
\rho_{k,k}(y)=\sum\limits_{n\geq
k+1}\sum\limits_{i=1}^{n-k}c_iop_{n-i,k}^{k+1}y^n=[C(y)-1]\varphi_{k+1}(y)=y^{k+2}[C(y)]^{k+4}.
\end{eqnarray*}

{\bf Case 2.} $k\geq m+1$
\begin{eqnarray*}\rho_{m,k}(y)&=&\sum\limits_{n\geq
k+1}op_{n,k}^my^n\\
&=&\sum\limits_{j=m+1}^{k+1}\sum\limits_{n\geq
k+1}\sum\limits_{i=1}^{n-k}c_iop_{n-i,k}^{j}y^n\\
&=&\sum\limits_{j=m+1}^{k}\sum\limits_{n\geq
k+1}\sum\limits_{i=1}^{n-k}c_iop_{n-i,k}^{j}y^n+\sum\limits_{n\geq
k+1}\sum\limits_{i=1}^{n-k}c_iop_{n-i,k}^{k+1}y^n\\
&=&[C(y)-1]\sum\limits_{j=m+1}^{k}\rho_{j,k}(y)+[C(y)-1]\varphi_{k+1}(y)\\
&=&y[C(y)]^2\sum\limits_{j=m+1}^{k}\rho_{j,k}(y)+y[C(y)]^2\varphi_{k+1}(y).\end{eqnarray*}

By the above equation, we have $\rho_{m+1,k}(y)
=y[C(y)]^2\sum\limits_{j=m+2}^{k}\rho_{j,k}(y)+y[C(y)]^2\varphi_{k+1}(y).$

Hence, $\rho_{m,k}(y)-\rho_{m+1,k}(y) =y[C(y)]^2\rho_{m+1,k}(y).$
Since $1+y[C(y)]^2=C(y)$, we obtain
$\rho_{m,k}(y)=C(y)\rho_{m+1,k}(y).$\hfill$\blacksquare$
\begin{corollary}
For any $m\geq 1$, let $\rho_{m}(x,y)$ be the generating function of
the sequence formed by $k$-flaw preference sets length $n$ with
leading term $m$. Then
\begin{eqnarray*}\rho_{m}(x,y)=\frac{x^my^{m+2}[C(y)]^{m+4}}{1-xy[C(y)]^2}.\end{eqnarray*}
Furthermore, let $\rho(x,y,z)$ be the generating function of the
sequence formed by $k$-flaw preference sets length $n$ with leading
term $m$. Then
\begin{eqnarray*}\rho(x,y,z)=\frac{xy^3z[C(y)]^5}{[1-xy[C(y)]^2][1-xyzC(y)]}.\end{eqnarray*}
\end{corollary}
{\bf Proof.} By Theorem \ref{generatingrecurrencemk}, we have
$\rho_{m,k}(y)=C(y)\rho_{m+1,k}(y)=[C(y)]^{k-m}\rho_{k,k}(y)=y^{k+2}[C(y)]^{2k-m+4}.$
Therefore, for any $k\geq m$, we have
$\rho_{m,k}(y)=y^{k+2}[C(y)]^{2k-m+4}.$

So,
\begin{eqnarray*}\rho_{m}(x,y)=\sum\limits_{k\geq m}y^{k+2}[C(y)]^{2k-m+4}x^k
=\frac{x^my^{m+2}[C(y)]^{m+4}}{1-xy[C(y)]^2}\end{eqnarray*} and
\begin{eqnarray*}\rho(x,y,z)=\sum\limits_{m\geq
1}\frac{x^my^{m+2}[C(y)]^{m+4}}{1-xy[C(y)]^2}z^m
=\frac{xy^3z[C(y)]^5}{[1-xy[C(y)]^2][1-xyzC(y)]}.\end{eqnarray*}\hfill$\blacksquare$\\

Let $\phi(x,y)=\sum\limits_{n\geq
1}\sum\limits_{k=0}^{n-1}op_{n,k}x^ky^n=\sum\limits_{k\geq
0}\sum\limits_{n\geq k+1}op_{n,k}y^nx^k$ be the generating function
for the sequence $op_{n,k}$. Let $\phi_k(x)=\sum\limits_{n\geq
k+1}op_{n,k}y^n$ for any $k\geq 0$, then
$\phi(x,y)=\sum\limits_{k\geq 0}\phi_k(y)x^k$.
\begin{theorem}\label{recurrenceopk-flawgenerating}
Let $k\geq 0$ and  $\phi_{k}(y)$ be the generating function  of the
sequence formed by $k$-flaw preference sets length $n$. Then
\begin{eqnarray*}\phi_k(y)=y^{k+1}[C(y)]^{2(k+1)}.\end{eqnarray*}
\end{theorem}
{\bf Proof.} By Theorem \ref{kflawrecurrence}, we have
\begin{eqnarray*}\sum\limits_{n\geq
k+1}op_{n,k}y^n&=&\sum\limits_{n\geq
k+1}\sum\limits_{i=0}^{n-k-1}c_i[op_{n-i,k}-op_{n-i,k}^1]y^n\\
&=&\sum\limits_{n\geq
k+1}\sum\limits_{i=0}^{n-k-1}c_iop_{n-i,k}y^n-\sum\limits_{n\geq
k+1}\sum\limits_{i=0}^{n-k-1}c_iop_{n-i,k}^1y^n.
\end{eqnarray*}
Hence, $\phi_k(y)=C(y)\phi_k(y)-C(y)\rho_{1,k}(y)$. Note that
$\rho_{1,k}(y)=y^{k+2}[C(y)]^{2k+3}$ and $C(y)-1=y[C(y)]^2$.
Therefore, $\phi_k(y)=y^{k+1}[C(y)]^{2(k+1)}.$ \hfill$\blacksquare$
\begin{corollary}\label{generaingkflawcorollary}
Let $\phi(x,y)$ be the generating function  of the sequence formed
by $k$-flaw preference sets length $n$. Then
\begin{eqnarray*}\phi(x,y)=\frac{y[C(y)]^2}{1-xy[C(y)]^2}.\end{eqnarray*}
\end{corollary}
{\bf Proof.} By Theorem \ref{recurrenceopk-flawgenerating}, since
$\phi(x,y)=\sum\limits_{k\geq 0}\phi_k(y)x^k$, we have
\begin{eqnarray*}\phi(x,y)=\sum\limits_{k\geq
0}y^{k+1}[C(y)]^{2(k+1)}x^k=\frac{y[C(y)]^2}{1-xy[C(y)]^2}.\end{eqnarray*}
\hfill$\blacksquare$

\section{The mean and the various of the number of flaws in $k$-flaw preference sets}
In this section, we will compute the mean and the various of the
number of flaws of the sequences corresponding with the generating
functions in Section $5$.

\begin{corollary}
For any  $k\geq 0$ and $n\geq k+1$, let $op_{n,k}^{k+1}$ denote the
number of $k$-flaw preference sets of length $n$ with leading term
$k+1$. Let $mop_{n,k}^{k+1}$ and $vop_{n,k}^{k+1}$ denote the mean
and various of  the number of flaws of preference sets in the set
$\mathcal{OP}_{n,k}^{k+1}$. Then we have
\begin{eqnarray*}mop_{n,k}^{k+1}=\frac{5(n-1)}{3(n+3)}, vop_{n,k}^{k+1}=\frac{4(n-1)(2n+1)(4n+7)}{9(n+3)^2(n+4)}.\end{eqnarray*}
\end{corollary}
{\bf Proof.} By Corollary \ref{generatingkflawleadingk+1}, the
generating function of the sequences formed by $op_{n,k}^{k+1}$ is
$\varphi(x,y)=\frac{y[C(y)]^2}{1-xyC(y)}$. Note that
$\frac{1}{C(y)}=1-yC(y)$. Simple computations tell us
\begin{eqnarray*}\left.\frac{\partial\varphi(x,y)}{\partial x}\right
|_{x=1}=y^2[C(y)]^5=y^{-1}\rho_{1,1}(y)=\sum\limits_{n\geq
2}\frac{5}{n+2}{2n-2\choose{n-3}}y^{n-1}=\sum\limits_{n\geq
1}\frac{5}{n+3}{2n\choose{n-2}}y^{n},\end{eqnarray*}

\begin{eqnarray*}\left.\frac{\partial^2\varphi(x,y)}{\partial x^2}\right
|_{x=1}=2y^3[C(y)]^7=2y^{-1}\rho_{1,2}(y)=\sum\limits_{n\geq
3}\frac{14}{n+3}{2n-2\choose{n-4}}y^{n-1}=\sum\limits_{n\geq
2}\frac{14}{n+4}{2n\choose{n-3}}y^{n},\end{eqnarray*} and
\begin{eqnarray*}\varphi(1,y)=y[C(y)]^3=y^{-1}\varphi_2(y)
=\sum\limits_{n\geq
2}\frac{3}{n+1}{2n-2\choose{n-2}}y^{n-1}=\sum\limits_{n\geq
1}\frac{3}{n+2}{2n\choose{n-1}}y^{n}.\end{eqnarray*}

Hence,
\begin{eqnarray*}mop_{n,m-1}^m=\frac{[y^n]\left.\frac{\partial\varphi}{\partial x}\right
|_{x=1}}{[y^n]\varphi(1,y)}=\frac{5(n-1)}{3(n+3)}\end{eqnarray*}

and
\begin{eqnarray*}vop_{n,m-1}^m&=&\frac{[y^n]\left.\frac{\partial^2\varphi}{\partial x^2}\right
|_{x=1}}{[y^n]\varphi(1,y)}+mop_{n,m-1}^m-[mop_{n,m-1}^m]^2\\
&=&\frac{4(n-1)(2n+1)(4n+7)}{9(n+3)^2(n+4)}.\end{eqnarray*}\hfill$\blacksquare$

\begin{corollary}For any $k\geq 0$, $m\geq k$ and $n\geq k+1$, let
$op_{n,k}^m$ denote the number of $k$-flaw preference sets of length
$n$ with leading term $m$. Let $mop_{n,k}^{m}$ and $vop_{n,k}^{m}$
denote the mean and various of  the number of flaws of preference
sets in the set $\mathcal{OP}_{n,k}^{m}$.Let
\begin{eqnarray*}
t(n,m)&=&\sum\limits_{i=0}^{n-m-2}\frac{m(m+3){2i\choose{i}}{2n-m-2-2i\choose{n-m-2-i}}+4^i(m+4){2n-m-3-2i\choose{n-m-3-i}}}{n-i+1}
\end{eqnarray*} and \begin{eqnarray*}
r(n,m)&=&\frac{m(m-1)(m+4)}{2n-m}{2n-m\choose{n-m-2}}+\sum\limits_{i=0}^{n-m-2}\frac{m+5}{2(n-i)-m-1}{2(n-i)-m-1\choose{n-i-m-3}}{2i\choose{i}}\\
&&+\sum\limits_{i=0}^{n-m-2}\frac{m+4}{2(n-i-1)-m}{2(n-i-1)-m\choose{n-i-m-3}}4^i.
\end{eqnarray*} Then  for $n\geq m+2$, we have
$$mop_{n,k}^m=\frac{t(n,m)}{{2n-m-1\choose{n-m-2}}}$$  and
$$vop_{n,k}^m=\frac{r(n,m)}{{2n-m-1\choose{n-m-2}}}+\frac{t(n,m)}{{2n-m-1\choose{n-m-2}}}-\left[\frac{t(n,m)}{{2n-m-1\choose{n-m-2}}}\right]^2.$$
\end{corollary}
{\bf Proof.} By Corollary \ref{generatingrecurrencemk}, the
generating function of the sequence formed by $op_{n,k}^m$ is
$\rho_m(x,y)=\frac{x^my^{m+2}[C(y)]^{m+4}}{1-xy[C(y)]^2}$. Note that
$\frac{1}{C(y)}=1-yC(y)$.

\begin{eqnarray*}\left.\frac{\partial\rho_m(x,y)}{\partial x}\right
|_{x=1}&=&\frac{my^{m+2}[C(y)]^{m+4}}{1-y[C(y)]^2}+\frac{y^{m+3}[C(y)]^{m+6}}{[1-y[C(y)]^2]^2}\\
&=&\frac{my^{m+2}[C(y)]^{m+3}}{\sqrt{1-4y}}+\frac{y^{m+3}[C(y)]^{m+4}}{1-4y}\\
&=&\frac{m\varphi_{m+2}(y)}{\sqrt{1-4y}}+\frac{\varphi_{m+3}(y)}{1-4y}\\
&=&\sum\limits_{n\geq m+2}t(n,m)y^n\end{eqnarray*}

\begin{eqnarray*}
\left.\frac{\partial^2\rho_m(x,y)}{\partial x^2}\right
|_{x=1}&=&m(m-1)y^{m+2}[C(y)]^{m+4}+\frac{2(m-1)y^{m+3}[C(y)]^{m+6}}{1-y[C(y)]^2}+\frac{2y^{m+3}[C(y)]^{m+6}}{[1-y[C(y)]^2]^2}\\
&=&m(m-1)y^{m+2}[C(y)]^{m+4}+\frac{2(m-1)y^{m+3}[C(y)]^{m+5}}{\sqrt{1-4y}}+\frac{2y^{m+3}[C(y)]^{m+4}}{1-4y}\\
&=&m(m-1)\rho_{m,m}(y)+\frac{2(m-1)\rho_{m+1,m+1}(y)}{\sqrt{1-4y}}+\frac{2y\rho_{m,m}(y)}{1-4y}\\
&=&\sum\limits_{n\geq m+2}r(n,m)y^n
\end{eqnarray*}

Since $\rho_m(x,y)=\sum\limits_{n\geq m+1}\sum\limits_{k=
m}^{n-1}op_{n,k}^mx^ky^n$, by Lemma \ref{allestkflawmaxlleadingm},
we have
\begin{eqnarray*}[y^n]\rho_m(1,y)=\sum\limits_{k=m}^{n-1}op_{n,k}^m=op_{n,\geq m}^m={2n-m-1\choose{n-m-2}}.\end{eqnarray*}

Hence, for $n\geq m+2$, we have
$$mop_{n,k}^m=\frac{t(n,m)}{{2n-m-1\choose{n-m-2}}}$$  and
$$vop_{n,k}^m=\frac{r(n,m)}{{2n-m-1\choose{n-m-2}}}+\frac{t(n,m)}{{2n-m-1\choose{n-m-2}}}-\left[\frac{t(n,m)}{{2n-m-1\choose{n-m-2}}}\right]^2$$\hfill$\blacksquare$

\begin{theorem}
For any $n\geq 1$, let $op_{n,k}$ denote the number of $k$-flaw
preference sets of length $n$. Let $mop_{n,k}$ and $vop_{n,k}$
denote the mean and various of  the number of flaws of preference
sets in the set $\mathcal{OP}_{n,k}$. Then we have
\begin{eqnarray*}
mop_{n,k}=\frac{2^{2n-1}}{{2n\choose{n}}}-1,
vop_{n,k}=n-\left[\frac{2^{2n-1}}{{2n\choose{n}}}\right]^2.
\end{eqnarray*}
\end{theorem}
{\bf Proof.} By Corollary \ref{generaingkflawcorollary}, the
generating function of the sequence formed by $op_{n,k}$ is
$\phi(x,y)=\frac{y[C(y)]^2}{1-xy[C(y)]^2}$. By Lemma \ref{B(y)},
simple computations tell us
\begin{eqnarray*}
\left.\frac{\partial\phi(x,y)}{\partial x}\right|_{x=1}&=&\frac{1}{2}+\frac{y}{1-4y}-\frac{1}{2\sqrt{1-4y}}\\
&=&\sum\limits_{n\geq
1}\left[4^{n-1}-\frac{1}{2}{2n\choose{n}}\right]y^n
\end{eqnarray*}
and
\begin{eqnarray*}
\phi(1,y)&=&-\frac{1}{2}+\frac{1}{2\sqrt{1-4y}} \\
&=&\sum\limits_{n\geq 1}\frac{1}{2}{2n\choose{n}}y^n.
\end{eqnarray*}
Hence,
\begin{eqnarray*}
mop_{n,k}=\frac{[y^n]\left.\frac{\partial\phi(x,y)}{\partial
x}\right|_{x=1}}{[y^n]\phi(1,y)}=\frac{2^{2n-1}}{{2n\choose{n}}}-1.
\end{eqnarray*}

It is easy to obtain the following equation.
\begin{eqnarray*}
\left.\frac{\partial^2\phi(x,y)}{\partial
x^2}\right|_{x=1}&=&\frac{1}{\sqrt{1-4y}}+\frac{y}{\sqrt{(1-4y)^3}}-\frac{3y}{1-4y}-1\\
&=&\sum\limits_{n\geq 1}\left[\frac{1}{2}(n+2){2n\choose{n}}-3\cdot
2^{2n-2}\right]y^n.
\end{eqnarray*}
Hence,
\begin{eqnarray*}
vop_{n,k}=\frac{[y^n]\left.\frac{\partial^2\phi(x,y)}{\partial
x^2}\right|_{x=1}}{[y^n]\phi(1,y)}+mop_{n,k}-(mop_{n,k})^2=n-\left[\frac{2^{2n-1}}{{2n\choose{n}}}\right]^2.
\end{eqnarray*}\hfill$\blacksquare$
\section{Appendix}
For convenience to check the equations in the previous sections, by
computer search, for $n\leq 6$, we obtain the number of $k$-flaw
preference sets of length $n$ with leading term $m$ and list them in
Table $1$. Note that $op_{n,k}=0$ if $n\leq k$ and $op_{n,k}^m=0$ if
$m\geq k+2$. When $k=0$, $op_{n,0}$ is equals the $n$-th Catalan
number $c_n$.

$$
\begin{array}{|r|r|l|l|l|l|l|l|l|l|}
\hline
 op_{n,k}^m&m=1&2&3&4&5&6&7&8&op_{n,k}\\
 \hline
 (n,k)=(1,0)&1&0&&&&&&&1\\
 \hline
 (2,0)&2&0&&&&&&&2\\
 \hline
 (2,1)&0&1&0&&&&&&1\\
 \hline
 (3,0)&5&0&&&&&&&5\\
 \hline
 (3,1)&1&3&0&&&&&&4\\
 \hline
 (3,2)&0&0&1&0&&&&&1\\
 \hline
 (4,0)&14&0&&&&&&&14\\
 \hline
 (4,1)&5&9&0&&&&&&14\\
 \hline
 (4,2)&1&1&4&0&&&&&6\\
 \hline
 (4,3)&0&0&0&1&0&&&&1\\
 \hline
 (5,0)&42&0&&&&&&&42\\
 \hline
 (5,1)&20&28&0&&&&&&48\\
 \hline
 (5,2)&7&6&14&0&&&&&27\\
 \hline
 (5,3)&1&1&1&5&0&&&&8\\
 \hline
 (5,4)&0&0&0&0&1&0&&&1\\
 \hline
 (6,0)&132&0&&&&&&&132\\
 \hline
 (6,1)&75&90&0&&&&&&165\\
 \hline
 (6,2)&35&27&48&0&&&&&110\\
 \hline
 (6,3)&9&8&7&20&0&&&&44\\
 \hline
 (6,4)&1&1&1&1&6&0&&&10\\
 \hline
 (6,5)&0&0&0&0&0&1&0&&1\\
 \hline
\end{array}
$$
\begin{center}
Table.1. $op_{n,k}^m$ for any $n\leq 6$
\end{center}

\end{document}